\documentclass{amsart}
\usepackage{graphicx}
\graphicspath{{figures/}}
\usepackage{amsmath,amsthm}
\usepackage{amssymb, verbatim}
\usepackage{tikz,lipsum,lmodern}
\usepackage[most]{tcolorbox}
\usepackage{multicol,wrapfig}
\usepackage{centernot}
\usepackage{sidecap,graphicx,subcaption}
\usepackage{hyperref}
\usepackage{cite}
\usepackage{float}

\DeclareMathOperator{\R}{\mathbb{R}}

\renewcommand{\P}{\mathbb{P}}

\newcommand{\tconv}{\mathsf{tconv}}
\newcommand{\tcone}{\mathsf{tcone}}

\newtheorem{theorem}{Theorem}

\newtheorem{proposition}[theorem]{Proposition}
\theoremstyle{definition}

\title{The tropical geometry of causal inference for extremes}

\author{Ngoc M Tran}
\address{Department of Mathematics, The University of Texas at Austin}
\email{tran.mai.ngoc@utexas.edu}
\thanks{The author would like to thank Bernd Sturmfels and Lior Pachter, whose requests for a written version of her talk has led to this paper. This work is supported by NSF Grant DMS-2113468 and the NSF
IFML 2019844 award to the University of Texas at Austin.}
\begin{document}

\begin{abstract}
Extreme value statistics is the max analogue of classical statistics, while tropical geometry is the max analogue of classical geometry. In this paper, we review recent work where insights from tropical geometry were used to develop new, efficient learning algorithms with leading performance on benchmark datasets in extreme value statistics. We give intuition, backed by performances on benchmark datasets, for why and when causal inference for extremes should be employed over classical methods. Finally, we list some open problems at the intersection of causal inference, tropical geometry and deep learning. 
\end{abstract}

\maketitle
\section{Introduction}
Extreme value statistics concerns the maxima of random variables and relations between the tails of distributions rather than averages and correlations. 
These are fundamental problems in many applications, from finance, engineering risks, nutrition to hydrology, to name a few \cite{janssen2020k, chautru2015, einmahl2018continuous, kluppelberg2021estimating, buck2021recursive,cooley2019decompositions}. 
Unique challenges to extreme values are lack of data and lack of model smoothness. As a demonstration, consider the extreme values analogue of the classical linear model $Y = BX$, the \emph{max-linear factor model} \cite{wang2011conditional,einmahl2012m}
\begin{equation}\label{eqn:max.linear.bigvee}
X_i = \max_{j=1,\dots,p}c_{ij}Z_j, \quad, i \in [d], c_{ij}, Z_j \geq 0, Z_j \mbox{ independent, unobserved},
\end{equation}
and its special case, the \textbf{max-linear Bayesian network} \cite{gissibl2018max,kluppelberg2019bayesian}
\begin{align}
X_i & =\max_{j=1,\dots,d}c_{ij}X_j \vee Z_j, \, i \in [d], c_{ij},  \label{eqn:max.linear.dag.bigvee} \\
s.t. & Z_j \geq 0, Z_j \mbox{ independent, unobserved}, \, C \mbox{ supported on a DAG } \mathcal{D}. \nonumber
\end{align}
These extreme value models are unavoidably nonlinear and not smooth, due to the presence of the max operator. Existing methods in extreme value statistics for dimensions $d \geq 3$ are limited due to an intractable likelihood, while existing techniques for learning Bayesian networks require a large amount of data to fit nonlinear models \cite{yu2019dag,zheng2018dags,zheng2020learning}. Therefore, models of \eqref{eqn:max.linear.bigvee} and \eqref{eqn:max.linear.dag.bigvee} require new learning methods.

A number of recent work have shown that certain statistical problems on extremes \emph{can be done well}. We discuss three specific examples in this paper.
\begin{enumerate}
  \item The Conditional Sampling Problem asks to draw samples from the distribution of $X_{\bar{K}}$ conditioned on $\{X_K = x_K\}$, $K \subset [d]$. Conditional sampling from extreme values are very difficult due to an intractable likelihood. Wang and Stoev \cite{wang2011conditional} gave an exact algorithm that allows the predictions of a learned max-linear model to be made at tens of thousands values on a personal laptop, a scale far out of reach for most other classes of models in extreme value statistics at that time \cite{davison2012statistical,davison2015statistics}. Their work has inspired similar sampling algorithms for parametrized models in extreme spatial statistics \cite{dombry2013regular} and time-series \cite{oesting2014conditional}. 
  \item The Latent River Problem asks to recover a river network from \emph{only extreme flows} from a set of $d$ stations, \emph{without} any information on the stations or terrain. Following the work of \cite{doi:10.1111/rssb.12355}, this became the benchmark problem for causal inference for extremes \cite{doi:10.1111/rssb.12355, mhalla2020causal}. The QTree method of Buck, Kl\"uppelberg and the author \cite{tran2021estimating} outperforms alternatives on multiple hydrology datasets, under the additional assumption that the causal network is a tree. 
  \item In \cite{amendola2022conditional}, Am\'{e}ndola, Lauritzen, Kl\"uppelberg and the author completely characterized the conditional independence theory for max-linear Bayesian networks. In particular, the work shows that this theory is \emph{different} from the traditional $d$-separation framework. Turning this theory into new algorithms for conditional independence testing for extreme values is an active research direction \cite{amendola2021markov}. 
\end{enumerate}

There is an opportunity to generalize these results to a complete theory of causal learning and inference for extremes. In this paper, we aim to convince the reader that tropical geometry may be an important ingredient for building such a theory. 
The Latent River problem in hydrology nicely illustrates all the core ideas of causal inference for extremes, so we will start there, in Section \ref{sec:case.study}. We explain why the max-linear Bayesian network is a natural extension of the classical Bayesian network for extremes, why this model is hard to learn. Contrary to the commonly held belief that extreme values model cannot be learned with classical methods, we demonstrate that very simple correlation-based methods can perform well on `nice' datasets such as the Danube. This calls for more rigorous selection of benchmark datasets for extreme values theory. In Section \ref{sec:qtree}, we discuss the tropical geometry and intuition behind QTree. In Section \ref{sec:background}, we briefly address two key questions: why max-linear model is important to extreme value statistics in general, and why tropical geometry naturally arise in their analysis. In Section \ref{sec:two.more}, we explain the tropical geometry behind the other aforementioned work, the Wang-Stoev algorithm \cite{wang2011conditional} and conditional independence theory for max-linear Bayesian networks \cite{amendola2022conditional}. Section \ref{sec:summary} concludes with some open directions. The author hopes that this perspective will open up a new interdisciplinary field between extreme values statistics and tropical geometry. 

\section{Case study: the Latent River problem}\label{sec:case.study}
Consider the problem of \emph{pollutant tracing} in subsurface waterways \cite{mcgrane2016impacts,WOLF20128} using sensor data \cite{bartos2018open,mao2019low}. That is, suppose we have sensors at $d$ locations that output observations $X^1, \dots, X^n \in \mathbb{R}_{\geq 0}^d$, where $X^t_i$ is the detected level of some pollutant at station $i$ at time $t$. We believe that the pollutants are transported by an unseen complex network of leaky sewage pipes and aquifers. Our goal is to estimate this network, so that when we observe large concentration of pollutants at a subset of the sensors, we can trace the source of that pollution event. The \emph{Latent River} problem is a benchmark version of this problem, where we need to learn a hidden river network given extreme water flow rates at a set of stations. Proposed methods are evaluated on how well the estimated river network approximate the true river network. Benchmark datasets used for causal inference for extermes so far included segments of the Danube river \cite{asadi2015extremes} and the Lower Colorado river \cite{tran2021estimating}. 

\begin{figure}[H]
\includegraphics[width=0.95\textwidth]{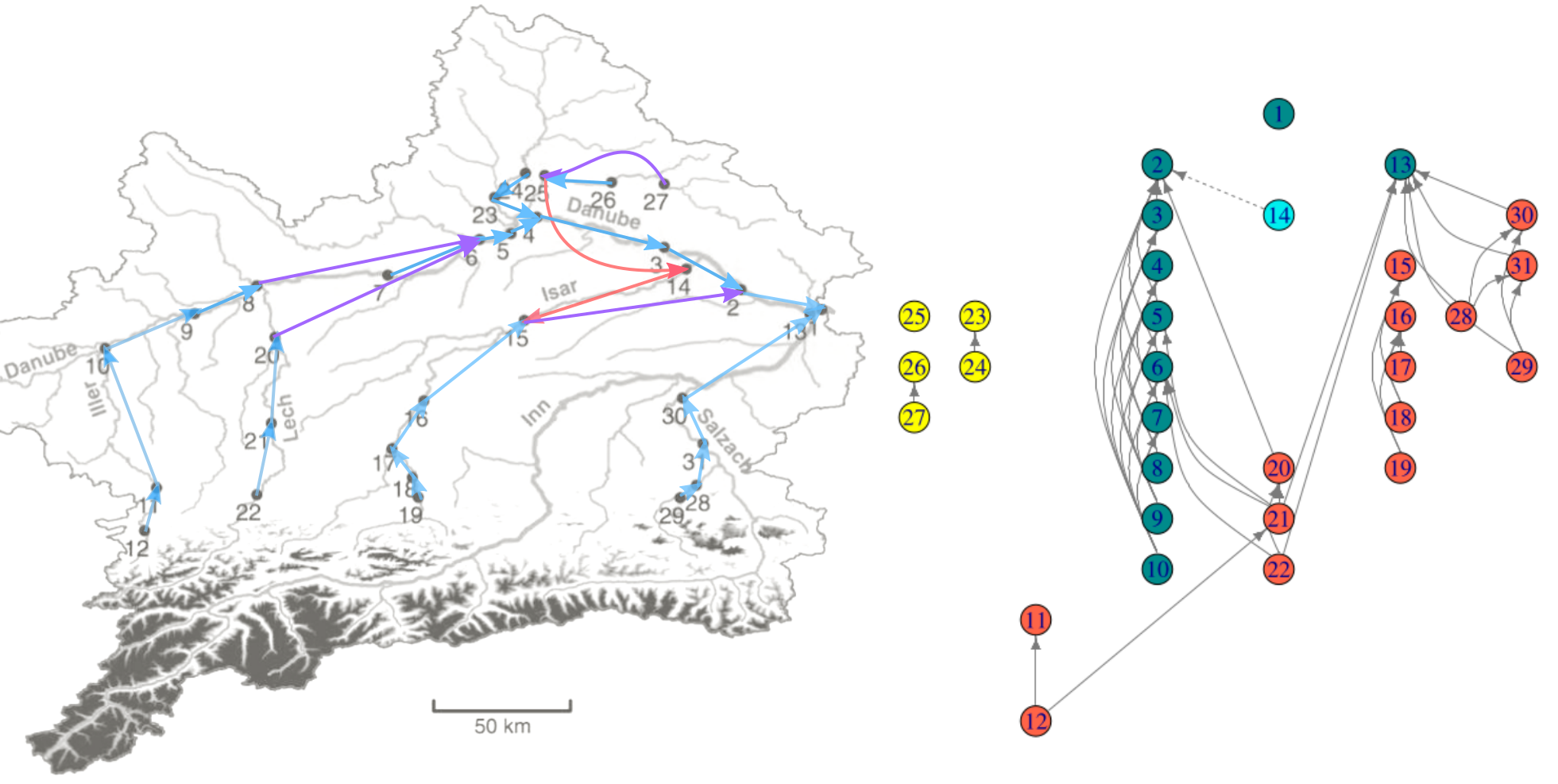}
\caption{The Hidden River Network asks to recover a river network from \emph{only extreme flows} from a set of $d$ stations, \emph{without} any information on the stations or terrain. This is a typical causal inference problem for extremes \cite{doi:10.1111/rssb.12355, mhalla2020causal}. \textbf{Left.} Map of the Danube (gray) \cite{asadi2015extremes} and the \emph{estimated, directed} river network obtained by QTree, a tropical algorithm for fitting a max-linear Bayesian tree \cite{tran2021estimating}. Correct edges are blue, wrong are red, purple edges are wrong edges with correct causal directions. \textbf{Right.} Performance of the previously leading method \cite{mhalla2020causal} on the same dataset, for comparison.} 
\label{fig:danube.with.competitor}
\end{figure}

Write $j \to i$ if node $j$ is directly upstream from node $i$, and $j \rightsquigarrow i$ if there is a directed path from $j$ to $i$. Let $G \in \{0,1\}^{d \times d}$ denote the unknown network, where $G_{ij} = 1$ iff $j \to i$. A simple Bayesian network model for this problem, ignoring time dependencies, would be
\begin{equation}\label{eqn:dag.linear}
X_i = \sum_{j: j \to i} c_{ij}X_j + Z_i + \epsilon_i, \quad c_{ij} > 0. 
\end{equation}
That is, pollutants at node $i$ consist of linear combinations of those from the parents nodes, plus those from an external source $Z_i$, plus measurement noise.

To make causal statements of the form `pollution at $i$ was caused by pollution at $j$` (and not the other way around), we need to assume that $G$ is a directed acyclic graph (DAG). Then \eqref{eqn:dag.linear} is a very standard Bayesian networks learning problem. While learning Bayesian networks is NP-complete \cite{chickering1996learning} since the space of DAGs is doubly exponential in the number of graph nodes, many algorithms have been developed. Score-based methods aim to optimize a score function such as BDeu and Bayesian information criterion over the space of DAGs using a search procedure \cite{ghahramani1997learning,murphy2002dynamic,chickering2002optimal,koivisto2004exact,silander2006simple,jaakkola2010learning,cussens2011bayesian,yuan2013learning,TALVITIE201969,niinimaki2016structure,pensar2020bayesian,kalisch2007estimating,bartlett2017integer,van2015machine}.
Recent work pursues a continuous optimization approach and grants more flexibility in both the choice of score function and the optimization procedure \cite{zheng2018dags,zheng2020learning}. 
When one of the inputs $Z_j$ (and consequently the measurement $X_j$) is large, the sum in \eqref{eqn:dag.linear} is dominated by the largest term, therefore \eqref{eqn:dag.linear} is approximately the max-linear Bayesian network \eqref{eqn:max.linear.dag.bigvee}
\begin{equation}\label{eqn:approx}
X = \sum_{j: j \rightsquigarrow i} c_{ij}X_j + Z_j \approx \bigvee_{j: j \rightsquigarrow i} c_{ij}X_j \vee Z_j + \mbox{ noise }.
\end{equation}
With this approximation, the Latent Tree is precisely the problem of learning a max-linear Bayesian network given noisy data. 

\subsection{Why learn from extremes. Why learning the Danube is easy but learning the Lower Colorado is hard.}

The problem of fitting the linear model \eqref{eqn:dag.linear} is unique to extreme values statistics in one of the following scenarios: 
\begin{enumerate}
	\item We only have access to extreme measurements in $X$. 
	\item We believe that there is a phase transition in the system, so that the network $G'$ that describes the causal relations between average measurements is fundamentally different from the network $G$ that describes the causal relations between extreme measurements. 
\end{enumerate}

For example, in hydrology, sensors can only reliable detect contaminants when the level is high enough. For benchmark data that uses water volume, rivers in drought such as the Lower Colorado can reigster $0$ or close to $0$ as the base flow at many stations. In this case, an edge $j \to i$ may be detectable \emph{only when} the flow is large enough to propagate from $j$ to $i$, and thus we are in scenario 2. Beyond hydrology, scenario 2 in particular is highly relevant in risk engineering \cite{gissibl2018max,gissibl2018tail,kluppelberg2019bayesian} and finance applications \cite{fama1993common,geluk2007weak,malevergne2002tail}. 

If the causal relations between extreme and typical values are the same, then one may as well fit the clasically linear Bayesian network \eqref{eqn:dag.linear} over \emph{typical} measurements instead of extreme values. By definition, typical measurements are more abundant than extreme ones. Furthermore, one can take advantage of the extensive body of knowledge developed for fitting classical linear Bayesian networks. In exchange, however, summation generally gives a distribution with higher entropy than the max, so a few extreme observations can carry more information on causal relations than a large number of averages. As a toy example, suppose we observe $X \in \R$, and we suspect that $X = Z_1 + Z_2$. Given $X = x$, we know very little about $Z_1$ and $Z_2$. Meanwhile, if $Y \in \R$ generated as $Y = Z_1 \vee Z_2$, then given $Y = y$, we know that \emph{either} $Z_1 = y$ or $Z_2 = y$, or both. Such exact inference techniques \cite{gissibl2018tail, gissibl2018max} do not work well for real data due to noise, however, one can generalize this observation to look for low-variance projections of the extreme values data that would give clues on the underlying DAG. This is the idea behind QTree \cite{tran2021estimating}, discussed in Section \ref{sec:qtree}. 

\begin{figure}[h]
\includegraphics[width=0.59\textwidth]{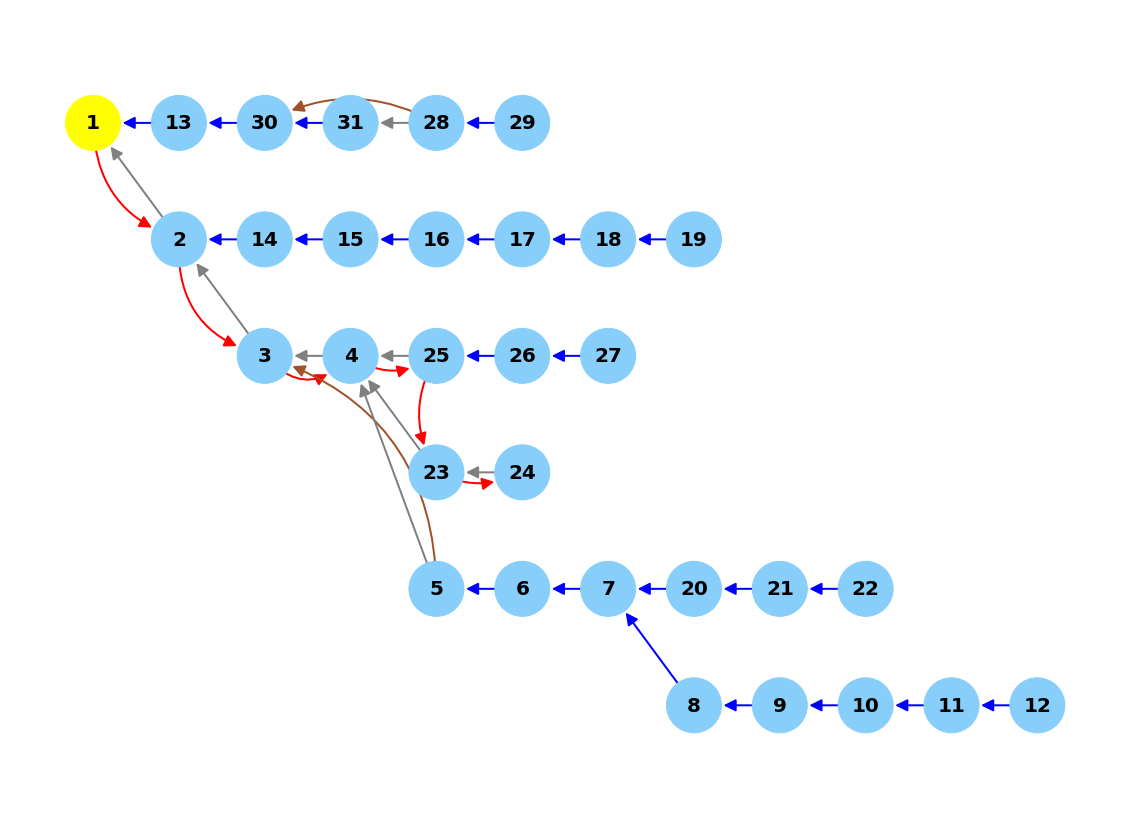}
\includegraphics[width=0.4\textwidth]{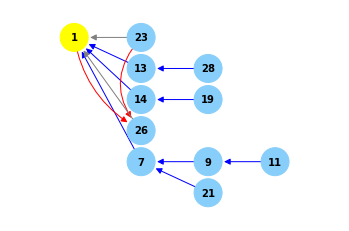} \\
\includegraphics[width=0.48\textwidth]{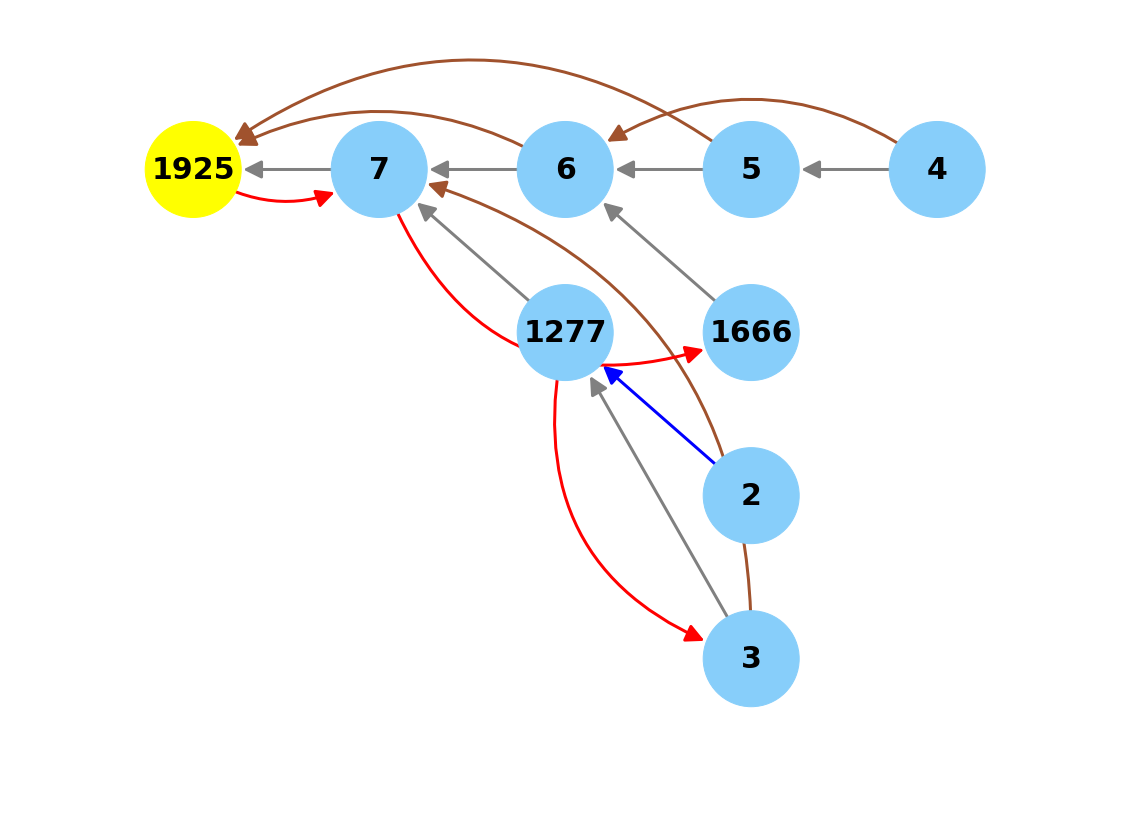}
\includegraphics[width=0.51\textwidth]{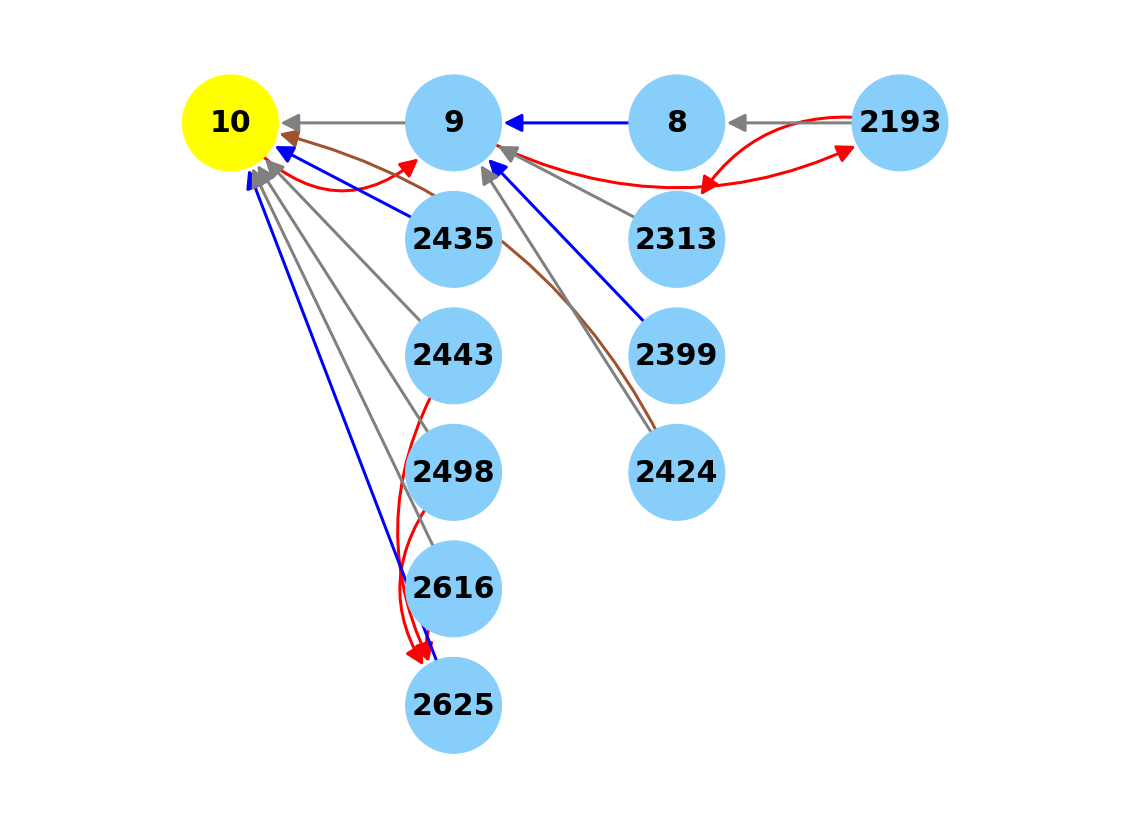} 
\caption{A na\"{i}ve correlation algorithm performs well on estimating the full Danube (top left) and an induced subset (top right), but fails for the upper (bottom left) and middle (bottom right) sections of the Lower Colorado. This is evidence that the Danube with smooth base flow is significantly easier to learn than the Lower Colorado, where base and extreme flows have different causal networks. Blue edges are correct, red edges are wrong, gray edges are correct edges absence from the estimated, and brown edges are wrong edges with correct causal direction. }\label{fig:danube_simple}
\end{figure}

Let us demonstrate the pros and cons of both methods on the Latent River problem. In extreme value statistics, the Danube river has emerged as the benchmark dataset for testing Latent Tree algorithms. We do not have access to the dataset of daily flows, only to the extreme measurements dataset created by Asadi, Davison and Engelke \cite{asadi2015extremes}. If we believe that the classical linear Bayesian network \eqref{eqn:dag.linear} is correct for both base flow and extreme flows, then these measurements, though extreme, are still drawn from this model. In the most na\"{i}ve algorithm, we simply compute the pairwise correlation of $X_i,X_j$, where $X_i$ is the observed flow at node $i$ and $X_j$ is the observed flow at node $j$. This gives a score matrix $S = (s_{ij})$, where higher score means it is more likely that $i$ and $j$ are neighbors on the true DAG $\mathcal{D}$. Then, under the assumption that the river is a tree, we run Edmond's algorithm to find the best-fitting directed tree that maximizes the score $S$. The result is shown in Figure \ref{fig:danube_simple}(top left). 

\begin{figure}[h]
\includegraphics[width=0.49\textwidth]{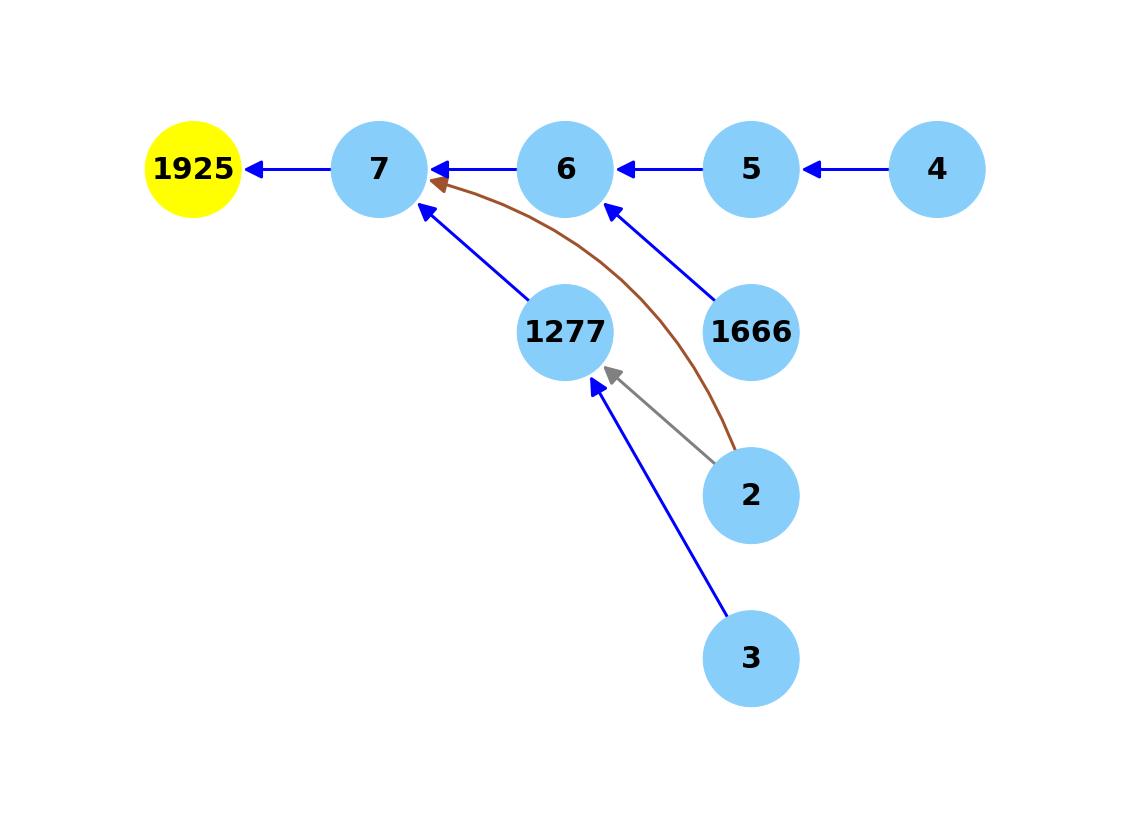}
\includegraphics[width=0.49\textwidth]{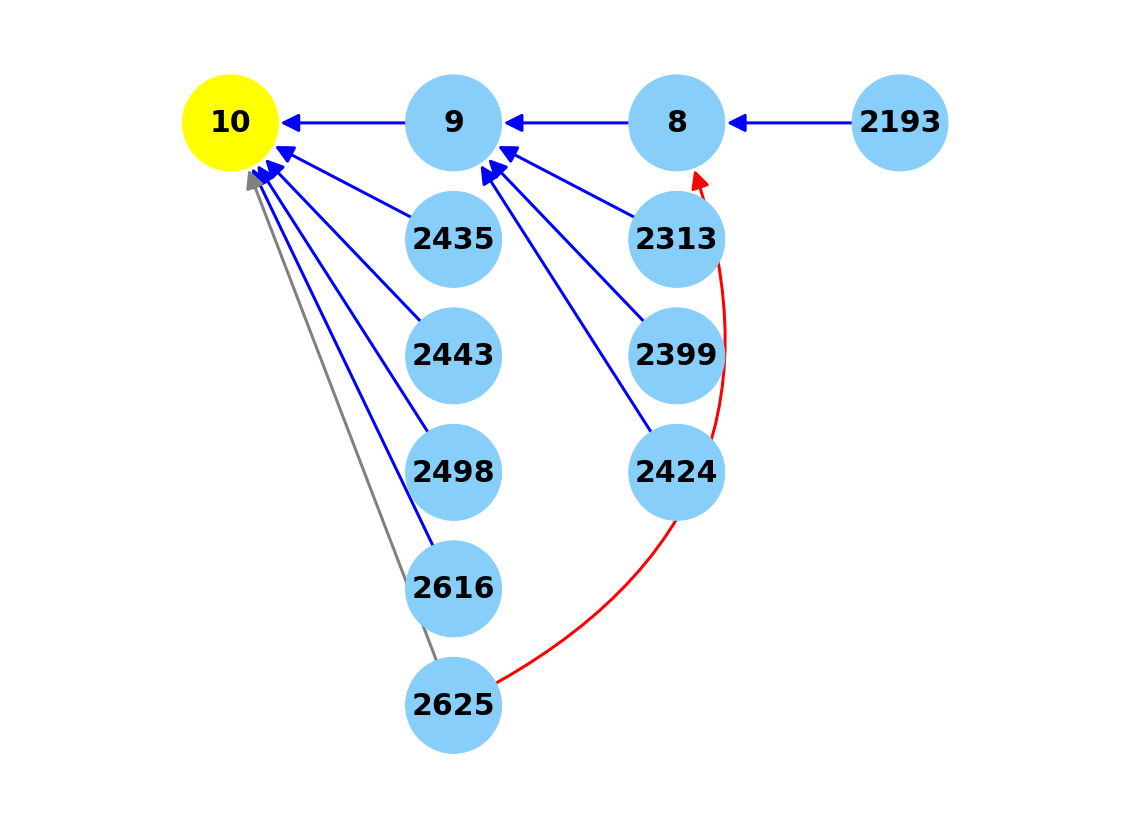} 
\caption{Upper (left) and Middle (right) sections of the Lower Colorado, estimated using QTree \cite{tran2021estimating} with default parameters.}\label{fig:qtree.colorado}
\end{figure}

On the full Danube, while this method gives more errors than the estimator returned by QTree, the result is comparable to those in the literature (cf. Figure \ref{fig:danube.with.competitor}(right)). In particular, it was able to correctly identify several sections of the Danube. We also apply the same algorithm on the same subset of nodes considered in \cite{gnecco2021causal}. The method of \cite{gnecco2021causal} gives perfect recovery on this subset and is undoubtedly more sophisticated. However, this simple algorithm yields a surprisingly good estimator, with only two incorrect edges. In comparison, the na\"{i}ve method fails to recover the Upper and Middle sections of the Lower Colorado. In comparison, QTree \cite{tran2021estimating} performs well on these datasets (cf. Figure \ref{fig:qtree.colorado}). Code and data necessary to produce these figures are at \url{https://github.com/princengoc/qtree/blob/main/qtree/naive_algorithm_vs_qtree.ipynb}. 

\begin{figure}[H]
\includegraphics[width=\textwidth]{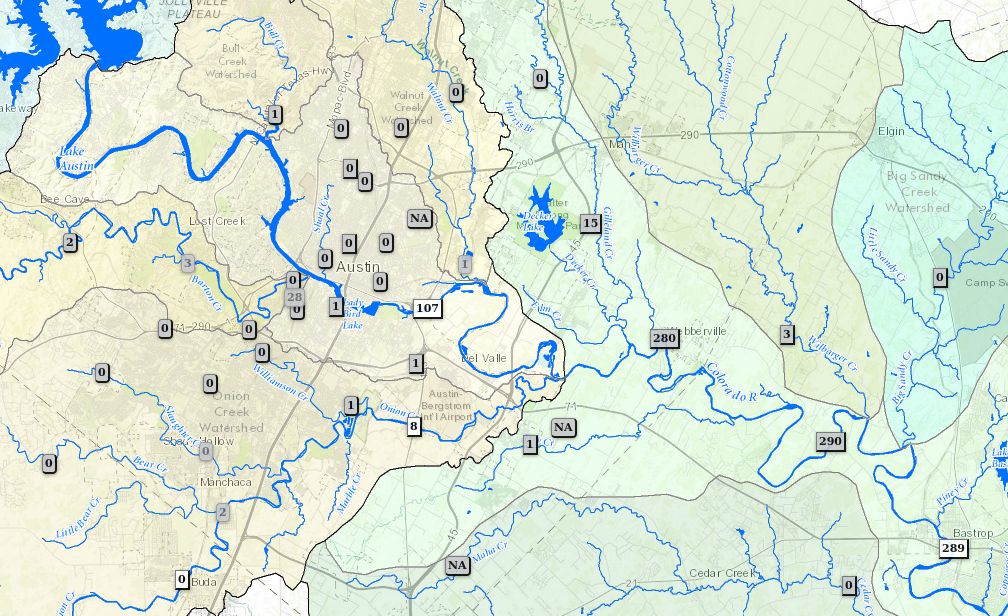} \\
\includegraphics[width=\textwidth]{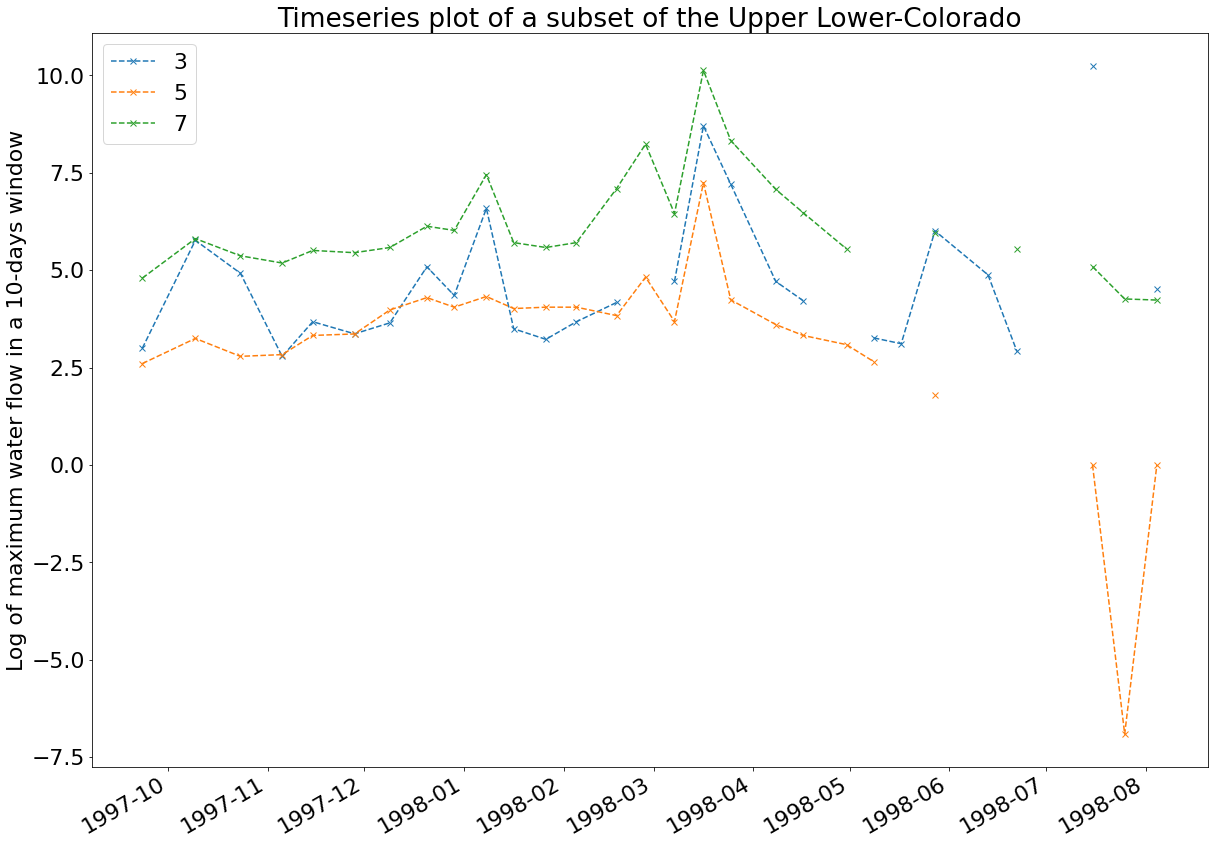}
\caption{\textbf{Top.} A subset of the Lower Colorado River and measurements at various stations, including minor stations off the main branch, snapshotted at a particular hour \cite{lcra}. \textbf{Bottom.} Three nodes on the dataset representing decorrelated extreme measurements in the Upper segment of the Lower Colorado \cite{tran2021estimating}, log-transformed. Zero measurements and missing data make causal discovery from extremes on the Lower Colorado much more challenging compared to the Danube.}\label{fig:lcra}
\end{figure}

The dataset of Upper, Middle and Lower Colorado are difficult. Each node has between 3\% and 55\% of data missing, largely caused by sensor failures \cite{lcra}. Some of these missing values are due to the sensor being hit by a flash flood, and thus the missing data point could occur exactly when an extreme measurement should have been recorded. But critically, unlike the Danube, these rivers are in perpetual drought and thus daily base flow can be close to 0 (cf. Figure \ref{fig:lcra}). As a result, two nearby stations may be essentially uncorrelated, and their link can \emph{only} be discovered when the upstream receives sufficiently large volume of water so that it can propagate downstream. For such dataset, causal inference methods built specifically for extreme values are needed. 
\newpage
\newpage

\subsection{Why learning max-linear Bayesian network is hard}

Section \ref{sec:background.extreme} discusses more extensively the difficulties in learning extreme values models in general. Here we briefly explain why `standard' methods for learning DAG on nonlinear models have not worked well for extreme values data such as the max-linear Bayesian network. First, we note that the vast majority of DAG fitting models have been developed for Gaussian or discrete data. General methods to learn DAG from nonlinear models with non-smooth objectives have only appeared in the last five years \cite{zheng2020learning, lachapelle2019gradient,yu2019dag}. However, they are designed for large neural networks and thus require a large amount of data and computation time. Extreme observations are few, and in applications such as hydrology, sensors can fail when a large volume of water hits, creating structurally missing data \cite{lcra} that is difficult to handle for gradient-based learning methods. However, we are optimistic that some modifications to the learning procedure of NOTEARS-DAG \cite{zheng2020learning} could work well. A different set of learning approaches is constraint-based methods. These aim to find the best DAG that respects the inferred conditional independence (CI) constraints
\cite{pearl2009causal,spirtes2000causation,zhang2008completeness}. Hybrid methods combine constraint-based and score-based approaches \cite{tsamardinos2003time,mooij2020joint}. However, these methods cannot be applied verbatim on max-linear Bayesian networks since they possess a \emph{different} CI theory that is \emph{not} characterized by $d$-separation \cite{amendola2022conditional}. In particular, computing CI statements for max-linear Bayesian networks from data is an active research area \cite{amendola2021markov}.  

Specific to the Latent River problem, a number of interesting methods have been proposed \cite{mhalla2020causal, engelke2021sparse, engelke2020graphical, gnecco2021causal}. On benchmark data that is Danube and Lower Colorado, to the best of our knowledge, the current leading method that can recover the whole directed Latent River is QTree \cite{tran2021estimating}, which heavily utilizes the \emph{tropical geometry} of the max-linear Bayesian network model. We discuss it next. 

\section{The tropical approach to Latent Tree}\label{sec:qtree}
Over the tropical \emph{max-times} semiring $(\R,\odot,\vee)$ where the max operator~$\vee$ replaces usual addition, the defining equation for the max-linear factor model \eqref{eqn:max.linear.bigvee} is a \emph{tropical matrix-vector multiplication}
\begin{equation}\label{eqn:max.linear.odot}
X = C \odot Z, \quad C = (c_{ij}) \in \R^{d \times p}_{\geq 0}, Z \in \R^p_{\geq 0}.
\end{equation}
The max-linear Bayesian network is a special case of \eqref{eqn:max.linear.odot} when $d = p$ and $C$ is a so-called Kleene star supported on a DAG \cite{tran2021estimating}. 
In geometric terms, \eqref{eqn:max.linear.odot} says that $X \in \R_{\geq 0}$ lies in the \emph{tropical cone} $\tcone(C)$ generated by $C$, which is the collection of all positive tropically linear combinations of the column vectors of the matrix $C$
$$\tcone(C) = \{C \odot z: z \in \R_{\geq 0}^p\}.$$
The tropical formulation is not only asthesically pleasing, it is \emph{far more efficient} to work with than the classical view. The log-normalization of $\tcone(C)$ is the tropical polytope $\tconv(C)$, which has~${d+p-2 \choose d-1}$ ordinary vertices \cite[Corollary 25]{develin2004tropical} as a piecewise-linear classical manifold, but only $p$ tropical vertices, which are the $p$ columns of $C$. In other words, the tropical algebra functions as a linearization tool, that turns the max-linear model from a piecewise linear, nonconvex object to a \emph{tropical polytope}. Learning \eqref{eqn:max.linear.odot} from noisy data is thus exactly the \emph{tropical nonnegative matrix factorization} (NMF). Geometrically, it is an approximate \emph{tropical convex hull} problem.

However, this rewriting alone does not lead to a straight-forward solution. Over both the classical and tropical algebra, determining the nonnegative rank of a matrix is NP-complete \cite{vavasis2010complexity, shitov2014complexity}. In other words, just learning the number if hidden factors $p$ in the max-linear model is hard. However, in some special cases, the max-linear model \emph{can} be learned well. The most obvious case is $d = p$, that is, $C = C^\ast$ is a tropical Kleene star. The max-linear Bayesian networks fall precisely into this category. In this case, idempotency of the tropical linear algebra leads to the result
$$ X = C^\ast \odot Z \mbox{ for some } Z \iff X = C^\ast \odot X. $$
Now, expanding the expression $X = C^\ast \odot X$ and take entry-wise $\log$, we find that if $j \rightsquigarrow i$, then
$$ \log(X)_i \geq \log(C)_{ij} + \log(X)_j, \iff \log(X)_i - \log(X)_j \geq \log(C)_{ij}. $$
That is, if $j$ is an ancestor of $i$, the log-measurements $\log(X)$, when projected onto the $e_i-e_j$ direction, forms a one-dimensional distribution with a \emph{lower bound} that is $\log(C)_{ij}$. Furthermore, equality happens whenever $X_i = C_{ij}X_j$, that is, whenever extreme measurement at $j$ and $i$ have a common cause. Thus, for noise-free data, we expect this one-dimensional distribution to have an atom at its minimum $\log(C)_{ij}$. With noisy data, we then expect an unusually large concentration of observations is some small quantile. The QTree algorithm \cite{tran2021estimating} crucially exploited this observation, using a local measure of variance as a measure of concentration and turn that into a score that reflects how confident we are that $j$ is an ancestor of $i$. More generally, the unique feature of data drawn from the noisy max-linear Bayesian network is that it has \emph{unusually low entropy in certain regions, under certain simple projections}. We are optimistic that this principle could be helpful in learning other piecewise linear models such as deep ReLU networks to noisy data, as deep ReLU networks are tropical rational functions \cite{zhang2018tropical}. 

\section{Why max-linear and why tropical}\label{sec:background}
Extreme value statistics is the max analogue of classical statistics, while tropical geometry is the max analogue of classical geometry. Thus, the geometry of extreme values statistics \emph{is} tropical. This section provides a short expansion of this statement and a brief review of each field's evolution. 

\subsection{Why Learning is difficult in extreme value statistics}\label{sec:background.extreme}
In probability theory max-stable distributions appear as limits of coordinate-wise maxima of i.i.d random variables, in the same way that Gaussians arise from sums of i.i.d's in the central limit theorem \cite{dehaan1970regular,de2007extreme}. Therefore, learning these distributions from data is a fundamental to extreme value statistics. While univariate and bivariate extreme value theory are fairly well understood and developed \cite{beirlant2006statistics, davison2015statistics}, learning and inference for max-stable multivariate in dimension $d \geq 3$ are {very difficult} due to an intractable likelihood. Indeed, for a max-stable random variable $X \in \R^d$ standardized to have unit Fr\'{e}chet marginals (an arbitrary choce), the CDF is \cite{beirlant2006statistics}
\begin{equation}\label{eqn:de.haan}
F(x) = \P(X_1 \leq x_1, \dots, X_d \leq x_d) = \exp(-\int_{\Delta_{d-1}}\max_{i=1,\dots,d}\frac{x_i}{w_i}S(dw)), x \in \R^d_{\geq 0},
\end{equation}
where $S$ is a positive measure on the simplex $\Delta_{d-1}$ known as the spectral measure of $X$. The likelihood is the derivative of \eqref{eqn:de.haan} with respect to $x_1,\dots,x_d$. This is a sum with $\sim d^{d-\log(d)}$ many terms, and thus is impractical to compute exactly in dimension $d \geq 3$ \cite{davison2012statistical}. This problem persists regardless of the parametrization of $S$, including the class of discrete spectral measures that correspond to the max-linear models \cite{fougeres2013dense}.

For $d \geq 3$, Bayesian and nonparametric methods exist but are very computationally expensive \cite{reich2012hierarchical,thibaud2016bayesian,padoan2022consistency}.
Most approaches with data applications follow one of two routes, as detailed in surveys \cite{davison2012statistical, davison2015statistics}. Composite likelihood approximates the likelihood with pairwise densities \cite{padoan2010likelihood}. Since max-linear models have discrete spectral measures, the likelihood is not smooth and thus causes difficulties \cite{einmahl2012m}. The other class of estimators aims to optimize the models' parameters directly through the empirical CDF of the original or transformed data \cite{segers2012nonparametric}. A typically transformed problem becomes that of estimating the \emph{tail dependence} function $\ell$ given $n$ approximate observations. However, the tail dependence function of a max-linear model is a sum of weighted maxima \cite{einmahl2012m} and thus is still not smooth. 
For max-linear models, one has \cite{einmahl2012m}
\begin{equation}\label{eqn:ell}
\ell(y) = \sum_{j=1}^p\max_{i=1,\dots,d}c_{ji}y_i.
\end{equation} 
Even with this transformation, the lack of smoothness of $\ell$ creates difficulties for likelihood approaches. Furthermore, some estimators based on the empirical CDF tend to suffer from numerical issues that worsen with $d$ and $p$, due to the need to numerically evaluate integrals over $[0,1]^d$ through discretization \cite{einmahl2012m,yuen2014crps}. Others sacrifice accuracy for speed, by focusing only on observations $y$ with large norms \cite{janssen2020k}, relying essentially on approximations for $\ell$ such as
$$ \ell(y) = \sum_{j=1}^p\max_{i=1,\dots,d}c_{ji}y_i \approx \max_{j=1,\dots,p}\max_{i=1,\dots,d}c_{ji}y_i $$
valid for $\|y\|$ extremely large. 
Recent literature has combined tail estimators with classical dimension reduction techniques in the form of principal components analysis \cite{haug2015copula, cooley2019decompositions,kluppelberg2021estimating}, clustering \cite{janssen2020k}, or multiple projections onto lower-dimensional cones under sparsity assumptions \cite{goix2017sparse,chiapino2019identifying}.  

\subsection{What is tropical geometry}\label{sec:tropical}
Tropical convex geometry and its algebraic counter part, tropical linear algebra, is geometry and algebra done over an idempotent semiring $(\R, \odot,\oplus)$, characterized by the property that $a \oplus a = a$ for all $a \in \R$. The two semirings relevant for this paper are the max-times, where $a \odot b := ab$, $a \oplus b := \max(a,b)$, and its log-analog, the max-plus, where $a \odot b := a+b$, $a\oplus b := \max(a,b)$. Originated from semigroup and finite automata theory \cite{mandel1977finite,simon1988recognizable,simon1994semigroups}, the field evolved through the 1980s as a linearization method in combinatorial optimization, with many applications to scheduling, queuing theory  \cite{baccelli1992synchronization, butkovivc2010max}, and large deviations \cite{puhalskii2001large}.
In the last five years, work by many authors has successfully used tropical linear algebra and convex geometry to solve important conjectures in classical linear programming \cite{allamigeon2018log} and economics \cite{baldwin2015understanding,tran2021finite}. Large-scale phylogenetics data \cite{yoshida2019tropical} along with the revelation that deep neural networks are tropical rational functions \cite{zhang2018tropical} have seen the emergence of the exciting field of \emph{tropical data science} \cite{yoshida2021tropical}.
Informally, the `tropical approach' to a problem involves two steps.
\begin{enumerate}
  \item Rewrite the problem (P) in a tropical algebra and reinterpret it as the classical problem (L) but with arithmetic done over the tropical algebra.  
  \item Apply the \emph{Maslov dequantization principle} \cite{litvinov1998linear}: if tool (T) is used to solve (L), use its tropical analog ($T_{trop}$) to solve (P).
\end{enumerate}

Step 1 is simply a change of dictionary. For example, learning the max-linear model becomes tropical nonnegative matrix factorization. The work starts with finding $T_{trop}$.  Loosely speaking, the Maslov principle says that one can continuously deform objects from the classical to the max-algebra such that key \emph{geometric} and \emph{combinatorial} features are preserved. The intuition of this step is that one can interpolate between the addition operator and the max operator, in the similar way that the $\ell_p$ norm family interpolates between $\ell_1$ and the $\ell_\infty$ \cite{litvinov1998linear}. Therefore, if tool (T) is a geometry-based solution to problem (L), then its tropical analog is likely to solve the original problem (P). Many theorems in tropical geometry formalize this intuition in specific settings \cite{maclagan2015introduction}. 

\subsection{How tropical geometry comes in to extreme value statistics}
Whenever one approximates a sum by its largest term, the max, then one has passed from classical to tropical geometry. This could happen as in our derivation of the max-linear Bayesian networks, which are arguably the simplest and most natural starting point for causal inference for extremes. 

Concretely, as of current, tropical geometry comes to extreme value statistics through the max-linear model (cf. Section \ref{sec:case.study}). While simple, this class of model is very important for two reasons. First, max-stable distributions can be approximated arbitrarily well by a max-linear factor model with i.i.d Frech\'{e}t sources $Z_j$'s \cite{stoev2005extremal, fougeres2013dense, falk2015}. Therefore, statistical theorems and techniques for max-linear models can be generalized to other max-stable distributions. These distributions, as aforementioned, are the Gaussians of extreme value statistics, their importance cannot be understated.   
Second, from a modeling perspective, max-linear models are the simplest class of models that exhibit \emph{cascading failure}, where extreme measurements (rainfall, contaminant level, risk, financial return) $X_i$ occur at a large number nodes can be traced to a few common sources (storm, chemical spill, catastrophic failure, financial shock) $Z_j$. Such cascading failures are commonly experienced in finance, engineering or hydrology, and therefore max-linear models are finding many applications in these domains \cite{einmahl2018continuous,gissibl2018graphical,gissibl2018max,tran2021estimating,kluppelberg2019bayesian,kluppelberg2021estimating, buck2021recursive, janssen2020k}. Indeed, max-linear models have long been used to model financial risks, starting with the extreme time-series max-ARMA model \cite{davis1989basic,davis1993prediction} and subsequent generalizations \cite{zhang2004behavior,ZHANG20062313,zhang2009approximating,ZHANG20101135,zhang2016copula,cui2018max}.
The max-linear factor model \eqref{eqn:max.linear.bigvee} is sometimes motivated as a way to study the \emph{tail dependence} function of the classical factor model $X = C Z + \epsilon$ popular in finance \cite{fama1993common,geluk2007weak,malevergne2002tail}, where $X$ is a multivariate regular variation, with heavy-tail factors $Z$ and light-weight noise $\epsilon$. 

\section{Two more examples of tropical geometry in extreme values theory}\label{sec:two.more}
\subsection{The Wang-Stoev conditional sampling algorithm}

For a random variable $X \in \R^d$, the Conditional Sampling Problem asks to draw samples from the distribution of $X_{\bar{K}}$ conditioned on $\{X_K = x_K\}$, $K \subset [d]$. Conditional sampling from extreme values distributions are very difficult due to an intractible likelihood. Wang and Stoev \cite{wang2011conditional} gave an exact algorithm that allows the \emph{predictions of a learned max-linear model} to be made at tens of thousands values on a personal laptop, a scale far out of reach for most other classes of models in extreme value statistics at that time \cite{davison2012statistical,davison2015statistics}. Their work has inspired similar sampling algorithms for parametrized models in extreme spatial statistics \cite{dombry2013regular} and time-series \cite{oesting2014conditional}. We now discuss the tropical geometry of ths algorithm. 

\begin{figure}[h]
\includegraphics[width=\textwidth]{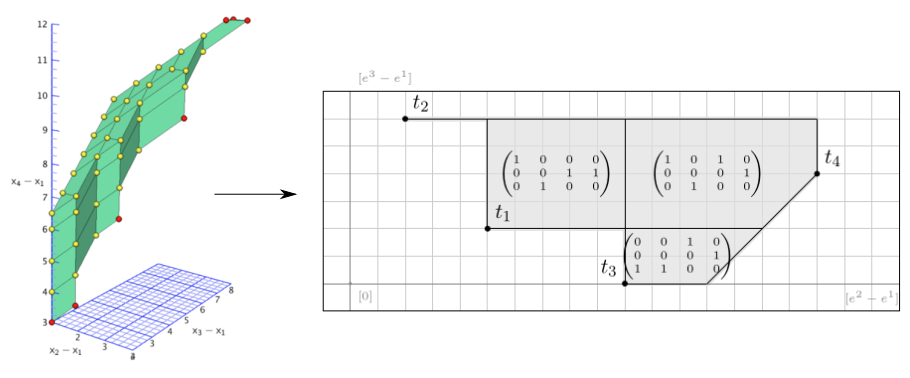}
\caption{The support of a max-linear factor model with $d = 4, p = 7$ (red points) is a tropical polytope $\tconv(C) \in \R^3$ (\textbf{left}, figure from \cite{Allamigeon2013MinimalER}). Conditioning on the event $\{X_K = x_K\}$ means slicing this figure with the subspace defined by $\{X_K = x_K\}$. The result is a lower-dimensional tropical polytope $\tconv(T)$ in dimension $\R^{d-|K|}$, in the variables $X_{\bar{K}} \in \R^{d-|K|}$. \textbf{Right}. An illustration of $\tconv(T)$ when $K = \{4\}$, with relatve covectors of its cells.}\label{fig:wang.stoev}
\end{figure}

Fix constants $C \in \R_{\geq 0}^{d \times p}$, $K \subset [d]$ and $x_K \in \R^{|K|}$. The Wang-Stoev algorithm is designed to sample from $X | X_K = x_K$, where $X = (X_{\bar{K}}, x_K) \in \R^d$ satisfies
$$ X = C \odot Z \quad \mbox{ $Z_j$'s independent and continuous, given } X_K = x_K. $$
As discussed in Section \ref{sec:qtree}, geometrically, the support of $X = C \odot Z \in \R^d$ is the tropical polytope $\tconv(C)$ generated by the columns of $C$. Then, $X_{\bar{K}} | X_K = x_K$ is obtained by slicing the original tropical polytope $\tconv(C)$ with an axis-aligned subspace (cf. Figure \ref{fig:wang.stoev}) and project onto the $\bar{K}$ coordinates. This is also a tropical polytope, which we will denote $\tconv(T)$.  Algebraically, $\tconv(T) \in \R^{\bar{K}}$ is the set of solutions to the tropically linear system
\begin{equation}\label{eqn:tropical.elimination}
\begin{pmatrix} X_{\bar{K}} \\ x_{K} \end{pmatrix} = 
\begin{pmatrix} 
C_{\bar{K}\bar{K}} & C_{\bar{K}K} \\ 
C_{K\bar{K}} & C_{KK}
\end{pmatrix} \odot 
\begin{pmatrix}
Z_{\bar{K}} \\ Z_{K}
\end{pmatrix}
\end{equation}
for known $C$ and $x_{K}$, and unknown $Z$. In other words, finding $\tconv(T)$, the support of the conditional distribution $X_{\bar{K}} | X_K = x_K$, is a \emph{tropical elimination} problem. 

Given $C$ and $x_K$, one can algorithmically solve for $\tconv(T)$ using, for example, the tropical Fourier-Motzkin algorithm \cite{allamigeon2014tropical}. The solution will directly express $X_{\bar{K}}$ as a nonlinear function of the $Z$'s and the constants $C$ and $x_K$, and from there one can draw samples from $X_{\bar{K}}$ given the joint distribution of $Z$. The advantage of this solution is generality: it works for \emph{any} joint distribution over the $Z_j$'s. The disadvatange is computational complexity. The tropical polytope $\tconv(T)$ is a union of \emph{cells} indexed by \emph{covectors}, which are certain bipartite matchings between $[d-K]$ and $[p]$ known in \cite{wang2011conditional} as \emph{hitting matrices}. Restricted to each cell, the map $Z \mapsto X_{\bar{K}}$ is linear, so it is always easy to sample from a cell of $\tconv(T)$ given its hitting matrix $H$. Thus one can imagine a two-step sampling algorithm: first compute a distribution $\nu_H$ over the hitting matrices. Then, draw a hitting matrix $H \sim \nu_H$, and from there, draw a sample $X_{\bar{K}}$ conditioned to be in the cell of $\tconv(T)$ with this hitting matrix. However, the number of valid hitting matrices grow as $O((d+p)^{d-|K|})$ for general $C$ \cite[Corollary 25]{develin2004tropical}. Wang and Stoev's main innovation \cite[Theorem 2]{wang2011conditional} is to show that when $Z_j$'s are independent, the probability distribution over hitting matrices $\nu_H$ has a nice structure that makes it easy to sample from, \emph{without} having to enumerate all the feasible hitting matrices. 

\subsection{Conditional independence theory of max-linear Bayesian networks}

Consider the max-linear Bayesian network  \eqref{eqn:max.linear.dag.bigvee} 
\begin{equation}\label{eqn:maxlinear2}
X = C \odot X \vee Z, \quad C \mbox{ supported on a DAG } \mathcal{D}.
\end{equation}
Suppose that $Z_j$'s are independent random variables. We want to answer questions of the form `is $X_I$ conditionally independent of $X_J$ given $X_K$`? Beyond intellectual curiosity, this question can lead directly to conditional independence tests, as well as algorithms for learning the DAG $\mathcal{D}$ given data, in the same way that classical conditional independence theory lead to causal inference algorithms of Gaussian models \cite{pearl2009causal,spirtes2000causation}. 

Write $X_I \perp X_J | X_K$ to denote the statement that $X_I$ is independent of $X_J$ given $X_K$. There are three levels at which we can answer this question. \\
\begin{enumerate}
  \item \emph{Context-specific}. Given $\mathcal{D}$, $C$ and $x_K$, is $X_I \perp X_J | X_K = x_K$? 
  \item \emph{Context-free, fixed $C$}. Given $\mathcal{D}$ and $C$, is $X_I \perp X_J | X_K = x_K$ true for \emph{all} possible values $x_K$? 
  \item \emph{Context-free, independent of $C$}. Given $\mathcal{D}$, is $X_I \perp X_J | X_K = x_K$ true for \emph{all} possible values $x_K$ and all possible $C$ supported on the DAG $\mathcal{D}$? 
\end{enumerate}

In \cite{amendola2022conditional}, Am\'{e}ndola, Kl\"{u}ppelberg, Lauritzen and the author give complete answers to these three problems, for all possible triples of sets of nodes $I,J,K \subseteq V$. Instead of reviewing the theorems or discuss their significance and implications to statistics, we will highlight the use of tropical geometry in this paper. 

The technical logic of \cite{amendola2022conditional} is as follows. Start with the context-specific case. We need to decide whether $X_I$ is independent of $X_J$, given $X_K = x_K$. The defining equation of the max-linear Bayesian networks \eqref{eqn:maxlinear2} has two equivalent forms:
$$ X = C \odot X \vee Z \iff X = C^\ast \odot Z \iff X = C^\ast \odot X, $$
where $C^\ast = I \vee \bigvee_{k=1}^d C^{\odot k}$ is the Kleene star matrix of $C$, and $I$ is the identity matrix. If we start with $X = C^\ast \odot Z$, then solving for $X_{\bar{K}}$ in terms of $Z_j$'s given $X_K = x_K$ is \emph{exactly} the tropical elimination problem we saw in \eqref{eqn:tropical.elimination}. Lucky for us, in this case, the Kleene star matrix $C^\ast$ has special structures that significantly simplify the geometry of the corresponding tropical polytopes $\tconv(C^\ast)$ and $\tconv(T)$, the support of $X_{\bar{K}} | X_K = x_K$. Algebraically, this geometry manifests in the fact that $X = C^\ast \odot Z$ for some $Z$ if and only if $X = C^\ast \odot X$. Using these equivalent forms, tropical elimination allows us to write the distribution of $X_{\bar{K}}$ as follows. 
\begin{proposition}[\cite{amendola2022conditional}, Proposition 4.1]\label{prop:main}
The following is a representation of $X_{\bar{K}} | X_K = x_K$ with respect to $Z$
$$ X_{\bar{K}} = C^\ast_{\bar{K}K} \odot x_K \vee C^\ast_{\bar{K}\bar{K}} \odot Z_{\bar{K}}, $$
where $Z$ is a multivariate distribution with independent coordinates conditioned to satisfy
\begin{equation}\label{eqn:prop.xk}
x_K = C^\ast_{KK}\odot Z_K \vee C^\ast_{K\bar{K}} \odot Z_{\bar{K}}.
\end{equation}
\end{proposition}
Intuitively, the $Z$'s are our sources of independence, but $x_K$ acts as a bridge that ties \emph{some} of them together, propagating `information` from one node to another, making them dependent. By definition, $C^\ast_{uv} > 0 \iff$ $v$ is an ancestor of $u$. Therefore, loosely speaking, any pair $Z_u,Z_v$ with at least one common descendant in $K$ will appear together in \eqref{eqn:prop.xk} and thus we would guess that they are dependent. Similarly, any pair $i,j \in K$ with a common ancestor in $\bar{K}$ will share a common term in the component $C^\ast_{K\bar{K}}Z_{\bar{K}}$, and thus we would again guess that they are dependent. These are `guesses` because we have not checked for potential cancellations due to the idempotency of the max operator. Indeed, the paper \cite{amendola2022conditional} offers several interesting examples, such as the Tent, the Umbrella and the Cassiopeia, where for exactly right values of $x_K$, cancellations occur that make seemingly dependent variables independent. Pathologies aside, the sketched intuition says that we should rewrite Proposition \ref{prop:main} as a block decomposition into common ancestors and descendants of $I$ and $J$, in and out of $K$ and $\bar{K}$. This resulted in \cite[Theorem 4.3]{amendola2022conditional}, and it is one of the main technical innovations of the paper. A more detailed decomposition featuring the tropical eigenvalue that led to sufficient and necessary conditions that cover all possible pathologies is given in \cite[Section 6]{amendola2022conditional}. 

\section{Summary}\label{sec:summary}
The geometry of max-linear models, the simplest of those in extreme value statistics, is tropical. In this paper, we review some recent work where tropical geometry, directly or indirectly, were used to achieve computationally efficient and accurate solutions to problems of learning and sampling from the max-linear model. We now collect some open research directions. 

\textbf{Problem 1. Find more benchmark datasets for causal inference for extremes}. The Danube dataset of Asadi, Davison and Engelke \cite{asadi2015extremes} have opened up an active research area. However, this work shows that the Danube may be a relatively simple dataset that cannot reliably distinguish between na\"{i}ve methods and sophisticated ones. Future benchmark datasets should ideally feature realistic challenges such as missing data, large number of nodes relative to the number of observations, a variety of DAG structures where `base flow` is different from `extreme flow`, and where causal statements made by a model could be independent tested and verified.  

\textbf{Problem 2. Generalize the Wang-Stoev conditional sampling algorithm}. A first step is to characterize all joint distributions over $Z$'s that give rise to computationally efficient $\nu_H$, the distribution over hitting matrices. This would lead to efficient conditional sampling algorithms for max-linear models with non-independent sources $Z_j$'s. It would also enable extensions to conditional sampling for other factor models in extreme values theory, parallel to the extension of the Wang-Stoev algorithm for parametrized models in extreme spatial statistics \cite{dombry2013regular} and time-series \cite{oesting2014conditional}. 

\textbf{Problem 3. Go beyond the max-linear model}. Whenever we approximate summation by maximum, we are passing from classical to tropical geometry. Such approximations are abundant in extreme value statistics. The max-linear model considered so far have given us tropical polytopes. However, both extreme value statistics and tropical geometry have powerful models and results beyond the linear case. What other possible collaborations could we discover between these fields? 

\bibliography{nsf.bib}
\bibliographystyle{plain}
\end{document}